\documentclass[11pt, reqno]{amsart}
\usepackage{amsmath, amssymb, amscd, amsthm, amsxtra}
\usepackage{graphicx}
\usepackage{color}
\usepackage{comment}

\usepackage[pdftex]{hyperref}

\usepackage[format=plain,labelfont=bf,up,width=.9\textwidth]{caption}

\theoremstyle{plain}
\newtheorem{thm}{Theorem}[section]
\newtheorem{lemma}[thm]{Lemma}
\newtheorem{prop}[thm]{Proposition}
\newtheorem{remark}[thm]{Remark}

\theoremstyle{definition}
\newtheorem{defn}[thm]{Definition}

\theoremstyle{plain}
\newtheorem{thmx}{\bf Theorem}

\newcommand{\R}{\mathbb{R}}
\newcommand{\N}{\mathbb{N}}
\newcommand{\C}{\mathbb{C}}
\newcommand{\Z}{\mathbb{Z}}
\newcommand{\D}{\mathbb{D}}
\newcommand{\Q}{\mathbb{Q}}
\newcommand{\s}{\mathbb{S}}
\newcommand{\bigO}{\mathcal{O}}

\newcommand{\Patt}{\mathcal{P}_{\rm att}}
\newcommand{\Prep}{\mathcal{P}_{\rm rep}}

\def\proof{\noindent {\bf Proof.\ }}

\def\qed{\hfill $\square$\\ }

\begin{document}
\title[Hedgehogs for neutral dissipative germs of $(\mathbb{C}^{2},0)$]{Hedgehogs for neutral dissipative germs of holomorphic diffeomorphisms of $(\mathbb{C}^{2},0)$}

\author[Tanya Firsova]{Tanya Firsova}
\address{Kansas State University, Kansas, United States}
\email{tanyaf@math.ksu.edu}

\author{Mikhail Lyubich}
\address{Stony Brook University, Stony Brook, United States}
\email{mlyubich@math.stonybrook.edu}

\author{Remus Radu}
\address{Stony Brook University, Stony Brook, United States}
\email{rradu@math.stonybrook.edu}

\author{Raluca Tanase}
\address{Stony Brook University, Stony Brook, United States}
\email{rtanase@math.stonybrook.edu}



\begin{abstract}
We prove the existence of hedgehogs for germs of complex analytic diffeomorphisms of $(\mathbb{C}^{2},0)$ with a semi-neutral fixed point at the origin, using topological techniques. This approach also provides an alternative proof of a theorem of P\'erez-Marco on the existence of hedgehogs for germs of univalent holomorphic maps of $(\mathbb{C},0)$ with a neutral fixed point.
\end{abstract}

\maketitle

\section{Introduction}\label{sec:Intro}

Let $\alpha\in\R\setminus\Q$ and let $p_{n}/q_{n}$ be the convergents of $\alpha$ given by the continued fraction algorithm. 
We say that $\alpha$ satisfies the Brjuno condition if 
\begin{equation}\label{eq:Brjuno}
\sum\limits_{n\geq 0}\frac{\log q_{n+1}}{q_{n}} < \infty.
\end{equation}
Brjuno \cite{Brj} and R\"{u}ssmann \cite{Rus} showed that if $\alpha$ satisfies Bjruno's condition, then 
any holomorphic germ with a fixed point with indifferent multiplier $\lambda = e^{2\pi i \alpha}$ is linearizable. The linearization is the irrational rotation with rotation number $\alpha$.
Yoccoz \cite{Y1} proved that Brjuno's condition is the optimal arithmetic condition that guarantees linearizability.  If $\alpha$ does not verify inequality \eqref{eq:Brjuno}, then there exists a holomorphic germ $f(z)=\lambda z+\bigO(z^{2})$ which is non-linearizable around the origin, that is $f$ is not conjugate to the linear map $z\mapsto \lambda z$ via a holomorphic change of coordinates. The origin is called a {\it Cremer} fixed point.

The local dynamics of a non-linearizable map with a Cremer fixed point is complex and hard to visualize. In the '90s, P\'erez-Marco \cite{PM1} proved the existence of interesting invariant compact sets near the Cremer fixed point, called hedgehogs. Using deep results from the theory of analytic circle diffeomorphisms developed by Yoccoz \cite{Y2}, P\'erez-Marco \cite{PM2, PM3, PM4} showed that  even if the map on a neighborhood of the origin is not conjugate to an irrational rotation, the points of the hedgehog are recurrent and still move under the influence of the rotation. 
Inou and Shishikura \cite{IS} built some models for the local dynamics near Cremer points for
specific cases of quadratic polynomials with high type rotation numbers using near-parabolic renormalization. 

In this paper we show the existence of non-trivial compact invariant sets for germs of diffeomorphisms of $(\C^2,0)$ with semi-indifferent fixed points.  The proof is purely topological and also provides an alternative proof for the existence of hedgehogs in dimension one.

A fixed point $x$ of a holomorphic germ $f$ of $(\C^{2},0)$ is {\it semi-indifferent} (or {\it semi-neutral}) if the eigenvalues $\lambda$ and $\mu$ of the linear part of $f$ at $x$ satisfy $|\lambda|=1$ and $|\mu|<1$. In analogy with the one-dimensional dynamics, a  semi-indifferent fixed point can be semi-parabolic, semi-Siegel or semi-Cremer, according to the arithmetic properties of the neutral eigenvalue $\lambda$.  We say that an isolated fixed point $x$ is {\it semi-parabolic}  if $\lambda=e^{2\pi i \alpha}$ and the angle $\alpha=p/q$ is rational. 
If $\alpha$ is irrational and there exists an injective holomorphic map $\varphi:\D\rightarrow \C^2$ such that $f(\varphi(\xi))=\varphi(\lambda \xi)$, for $\xi\in\D$, we call the fixed point {\it semi-Siegel}. Finally, if $\alpha$ is irrational and there does not exist an invariant disk on which the map is analytically conjugate to an irrational rotation, then the fixed point is called {\it semi-Cremer}. Note that in the latter case $\alpha$ does not satisfy the Brjuno condition \eqref{eq:Brjuno}.

Let $E^s$ and $E^c$ denote the eigenspaces of $Df_0$ corresponding to the dissipative eigenvalue $\mu$ and respectively to the neutral eigenvalue $\lambda$. Let $B$ be a neighborhood of $0$ and $E^s_x$ and $E^c_x$ not necessarily invariant continuous distributions such that $E^s_0=E^s$,  $E^c_0=E^c$, and $T_xB=E^s_x\oplus E^c_x$ for all $x\in B$. We define the vertical cone $\mathcal{C}^{v}_{x}$ to be the set of vectors in the tangent space at $x$ that make an angle less than or equal to $\alpha$ with $E^{s}_x$, for some $\alpha>0$.  The horizontal cone $\mathcal{C}^{h}_{x}$ is defined in the same way, with respect to $E^{c}_x$. 

The map $f$ is {\it partially hyperbolic} on $B$ (see Pesin \cite{P}) if there exist two real numbers $\overline{\mu}$ and $\underline{\lambda}$ such that $0<|\mu|<\overline{\mu}<\underline{\lambda}<1$ and a family of invariant cone fields $\mathcal{C}^{h/v}$ on $B$, 
\begin{equation}\label{eq:inv}
Df_x(\mathcal{C}_x^{h})\subset \mbox{ Int } \mathcal{C}_{f(x)}^{h}\cup\{0\}, \ \ \ Df^{-1}_{f(x)}(\mathcal{C}_{f(x)}^{v})\subset \mbox{ Int } \mathcal{C}_{x}^{v}\cup\{0\},
\end{equation}
such that for every $x\in B$ we have
\begin{equation}\label{eq:Cn}
 \underline{\lambda}\,\|v\|\leq \|Df_x(v)\|\leq 1/\underline{\lambda}\,\|v\|,\ \ \ \ \ \mbox{for}\ v\in \mathcal{C}_x^{h},
 \end{equation}
\begin{equation}\label{eq:Cs}
  \|Df_x(v)\|\leq \overline{\mu}\, \|v\|,\ \ \ \ \ \mbox{for}\ v\in \mathcal{C}_{x}^{v},
\end{equation}
for some Riemannian metric $\|\cdot\|$. 

If $f$ is partially hyperbolic on $B$, then the rate of contraction along $E^s_x$ dominates the behavior of $Df_x$ along the complementary direction $E^c_x$. This domination ensures the existence of 
local center manifolds $W^c_{\rm loc}(0)$ relative to $B$ as graphs of functions $\varphi_f: E^c\cap B\rightarrow E^s$, as discussed in Section \ref{sec:centermfd}.

\begin{thmx}\label{thm:Hedgehog} 
Let $f$ be a germ of holomorphic diffeomorphisms of $(\C^{2},0)$ with a semi-indifferent fixed point at $0$ with eigenvalues $\lambda$ and $\mu$, where $|\lambda|=1$ and $|\mu|<1$. Consider an open ball $B\subset \C^{2}$ centered at $0$ such that $f$ is partially hyperbolic on a neighborhood $B'$ of $\overline{B}$.

There exists a set $\mathcal{H}\subset \overline{B}$ such that:
\begin{itemize}
\item[a)] $\mathcal{H}\Subset W^{c}_{\rm loc}(0)$, where $W^{c}_{\rm loc}(0)$ is any local center manifold of the fixed point $0$ corresponding to the neutral eigenvalue $\lambda$, constructed relative to $B'$.
\item[b)] $\mathcal{H}$ is compact, connected, completely invariant and full.
\item[c)] $0\in \mathcal{H}$, $\mathcal{H}\cap\partial B\neq \emptyset$. 
\item[d)] Every point $x\in\mathcal{H}$ has a well defined local strong stable manifold $W^{ss}_{\rm loc}(x)$, consisting of points from $B$ that converge asymptotically exponentially fast to $x$.
The strong stable set of $\mathcal{H}$ is laminated by vertical-like holomorphic disks.
\end{itemize}
\end{thmx}

We say that $\mathcal{H}$ is {\it completely invariant} if $f(\mathcal{H})\subset \mathcal{H}$ and $f^{-1}(\mathcal{H})\subset \mathcal{H}$.  The set $\mathcal{H}$ is {\it full} if its complement in $W^{c}_{\rm loc}(0)$ is connected.  The local strong stable manifold $W^{ss}_{\rm loc}(x)$ of a point $x\in\mathcal{H}$ is defined as the set 
\begin{equation*}\label{eq:Wssx}\{y\in B : f^{n}(y)\in B\ \forall n\geq 1,\ \lim\limits_{n\rightarrow \infty} \mbox{dist}(f^n(y),f^n(x))/\overline{\mu}^{n} = 0\},
\end{equation*}
where $\overline{\mu}$ is the constant of partial hyperbolicity from \eqref{eq:Cs}.

We call the set $\mathcal{H}$ from Theorem \ref{thm:Hedgehog} a {\it hedgehog}. The most intriguing case is when  the argument of $\lambda$ is irrational and $\mathcal{H}$ is not contained in the closure of a linearization domain. This happens for instance, when the origin is semi-Cremer. Theorem \ref{thm:Hedgehog} is applicable to the local study of dissipative polynomial automorphisms of $\C^2$ with a semi-indifferent fixed point.

The theorem generalizes directly to the case of holomorphic germs of diffeomorphisms of $(\C^n,0)$, for $n>2$, which have a fixed point at the origin with exactly one eigenvalue on the unit circle and $n-1$ eigenvalues inside the unit disk. 

\medskip
\noindent\textit{Acknowledgement.} We are grateful to Romain Dujardin for reading an earlier version of the paper and for providing valuable comments. The first author was supported by NSF grant DMS-1505342. The second author was supported by NSF grants DMS-1301602 and DMS-1600519.

\section{Center manifolds of the semi-indifferent fixed point}\label{sec:centermfd}

Let $f:(\C^2,0)\rightarrow(\C^2,0)$, $f(x,y)=(\lambda x + f_1(x,y), \mu y + f_2(x,y))$ be a holomorphic germ with a semi-indifferent fixed point at the origin.
We also refer to $f$ as a neutral dissipative germ of $(\C^2,0)$. 

The semi-indifferent fixed point has a well-defined unique analytic strong stable manifold $W^{ss}(0)$ corresponding to the dissipative eigenvalue $\mu$. It consists of points that are attracted to $0$ exponentially fast, and defined as
\begin{equation}\label{eq:Wss}
W^{ss}(0):=\{x\in\C^{2} : \lim\limits_{n\rightarrow \infty} \mbox{dist}(f^{n}(x), 0)/\mu^{n}=\mbox{const}.\}.
\end{equation}

The semi-indifferent fixed point also has a (non-unique) center manifold $W^{c}_{\rm loc}(0)$ of class $C^k$ for some integer $k\geq1$, tangent at $0$ to the eigenspace $E^c$ of the neutral eigenvalue $\lambda$. There exists a ball $B_{\delta}$ (where the size of $\delta$ depends on $k$) centered at the origin in which the center manifold is locally the graph of a $C^k$ function $\varphi_{f}: E^c\rightarrow E^s$ and has the following properties:
\begin{itemize}
\item[a)] \textbf{Local Invariance:} $f(W^c_{\rm loc}(0))\cap B_{\delta}\subset W^c_{\rm loc}(0)$. 
\item[b)] \textbf{Weak Uniqueness:} 
If $f^{-n}(x)\in B_{\delta}$ for all $n\in\N$, then $x\in W^c_{\rm loc}(0)$. Thus center manifolds may differ only on trajectories that leave the neighborhood $B_{\delta}$ under backward iterations. 
\item[c)] \textbf{Shadowing:} Given any point $x$ such that $f^{n}(x)\rightarrow 0$ as $n\rightarrow \infty$, there exists a positive constant $k$ and a point $y\in W^c_{\rm loc}(0)$ such that $\|f^{n}(x)-f^{n}(y)\|<k\overline{\mu}^n$ as $n\rightarrow \infty$. In other words, every orbit which converges to the origin can be described as an exponentially small perturbation of some orbit on the center manifold.  
\end{itemize}

Consider the space of holomorphic germs $g$ of $(B_{\delta},0)$ which are $C^k$-close to $f$ such that $g$ has a semi-indifferent fixed point at the origin. We will later consider a sequence of germs with a semi-parabolic fixed point which converges uniformly to a germ with a semi-Cremer fixed point. Proposition \ref{prop:family} shows that even if the center manifold of $g$ is not unique, it may be chosen to depend continuously on $g$ for the $C^k$ topology. Let $E^{c}(\delta)=E^{c}\cap B_{\delta}$.

\begin{prop}\label{prop:family} The map $g$ has a $C^k$ center manifold defined as the graph of a $C^k$ function $\varphi_g: E^c(\delta)\rightarrow E^s(\delta)$ such that the map $(g,x)\mapsto \varphi_g(x)$ is $C^k$ with respect to $g$. 
\end{prop}

We refer to  \cite{S} (Chapter 5, Appendix III) and \cite{HPS} for the theory of stable and center manifolds.
For a proof of Proposition \ref{prop:family} see Theorem 5.1 and \S5A in \cite{HPS}. 

\smallskip
Assume that the map $f$ is partially hyperbolic on an open set $B'$, as in Theorem \ref{thm:Hedgehog}. 
We consider local center manifolds defined with respect to the ball  $B_{\delta}=B'$, satisfying the three properties above. 
By using a cut-off function we can construct a $C^{k}$-smooth extension $\tilde{f}$ of $f$ to $\C^{2}$
and make the center manifold $W^{c}(0)$ globally defined (see \cite{R}, \cite{V}). 
The extension $\tilde{f}$ can be chosen such that $\tilde{f}=f$ on $B_{\delta}$, $\tilde{f}=Df_0$ on the complement of $B_{2\delta}$ in $\C^{2}$, and $\|\tilde{f}-Df_{0}\|_{C^{1}}<\epsilon$ for some small $\epsilon$, which depends on the constants of partial hyperbolicity $\overline{\mu}$ and $\underline{\lambda}$ from  \eqref{eq:Cn} and \eqref{eq:Cs}.
 The proof of the existence of the center manifold for the modified function follows the usual contracting argument: consider the space of graphs of Lipschitz maps $h:E^c\rightarrow E^s$. The fact that the strong contraction along $E^s$ (and in the vertical cones $\mathcal{C}^v_x$, $x\in\C^{2}$) dominates the behavior of $Df$ along $E^c$ (and in the horizontal cones $\mathcal{C}^h_x$) ensures that the action of $\tilde{f}$ on the space of graphs is a contraction, hence it has a unique fixed point, which is the center manifold $W^c(0)$. This is globally defined and homeomorphic to $\R^2$, but clearly non-unique for the initial function $f$ as it depends on the choice of the extension function $\tilde{f}$. 
 
 To obtain the $C^k$-smoothness of the center manifold, it suffices to assume that the constants $\overline{\mu}$ and $\underline{\lambda}$ 
 satisfy $\overline{\mu} \underline{\lambda}^{-j}<1$ for $1\leq j \leq k$ on $B_{\delta}$, condition which is true if $\delta$ is small, since $|\mu|<1$ and $|\lambda|=1$. We will only use center manifolds of class $C^1$, so the condition $0<\overline{\mu}<\underline{\lambda}<1$ in the definition of partial hyperbolicity suffices for our purposes.

\section{Semi-parabolic germs}\label{sec:semiparabolic}
In this section we discuss the local structure of a germ $f$ of holomorphic diffeomorphisms of $(\C^2,0)$ with a semi-parabolic fixed point at the origin and we show the existence of big invariant petals. 
For simplicity, we call $f$ a {\it semi-parabolic germ}. 
Throughout the section, denote the eigenvalues of $Df_0$ by $\lambda=e^{2\pi i p/q}$ and $\mu$, where $|\mu|<1$ and $p/q$ is a rational number with $gcd(p,q)=1$. The following result is Proposition 3.3 from \cite{RT}. 

\begin{prop}\label{thm:normalform3}
Let $f$ be a semi-parabolic germ of transformation of $(\C^{2},0)$, with eigenvalues $\lambda$ and $\mu$, with $\lambda=e^{2\pi i p/q}$ and $|\mu|<1$. There exists a neighborhood $U$ of $0$ and local coordinates $(x,y)$ on $U$ in which $f$ has the form $f(x,y)=(x_{1},y_{1})$, where
\begin{equation}\label{eq:NF3}
\left\{\begin{array}{l}
    x_{1}= \lambda (x+x^{\nu q+1} + Cx^{2\nu q+1}+ a_{2\nu q+2}(y)x^{2\nu q +2}+\ldots )\\
    y_{1}= \mu y + xh(x,y)
\end{array}\right.
\end{equation}
and $C$ is a constant. The multiplicity of the fixed point as a solution of the equation $f^{q}(x,y)=(x,y)$ is $\nu q +1$. We call $\nu$ the semi-parabolic multiplicity of the fixed point.
\end{prop}

In this section, we will only deal with semi-parabolic multiplicity $\nu=1$. 
Otherwise we would have $\nu$ cycles of $q$ petals, invariant under $f^{q}$. Using the results of Ueda \cite{U1,U2} and Hakim \cite{Ha}, we can describe the local dynamics of semi-parabolic germs as follows:

\begin{thm}\label{thm:InvPetals} 
Let $f$ be a semi-parabolic germ of transformation of $(\C^{2},0)$, with eigenvalues $\lambda$ and $\mu$, with $\lambda=e^{2\pi i p/q}$ and $|\mu|<1$. Assume that the semi-parabolic multiplicity is $1$.
Let $U$ be the normalizing neighborhood from Proposition \ref{thm:normalform3}. Inside $U$, there exist $q$ attracting petals $\mathcal{P}_{{\rm att}, j}$ and $q$ repelling petals $\mathcal{P}_{{\rm rep}, j}$ for $1\leq j\leq q$. The attracting petals are two-dimensional and are a local base of convergence for $f$. The repelling petals are one-dimensional and are a local base of convergence for $f^{-1}$. 
\end{thm}

There are several ways to define local attractive and repelling petals.
 We will define fat attractive petals as in \cite[Section~4]{RT}. Consider the sets 
\begin{equation}\label{eq:Delta}
\Delta_{r}^{\pm} =\{x\in \C:\left({\rm Re}(x)\pm r\right)^{2}+\left(|{\rm Im}(x)|-r\right)^{2}<2r^{2}\}.
\end{equation}
and the coordinate transformation $T:(x,y)\mapsto (x^{q},y)$. Geometrically, $\Delta_r^+$ is the union of two complex disks of radius $\sqrt{2}r$ centered at $-r\pm ir$. Let $\mathcal{P}_{{\rm att}, j}$, $1\leq j\leq q$, be the preimages under $T$ of the set $\{x\in \Delta_{r}^{+}, |y|<r'\}$, for some $r,r'>0$ sufficiently small.  

A construction of repelling petals is done by Ueda in \cite{U2} for $\lambda=1$. Ueda shows that the local repelling petal is a smooth graph $\{y=\psi(x), |x-r|<r\}$, for some $r>0$ sufficiently small.
However, using the same techniques from \cite{U2} we can construct a larger local repelling petal which is a smooth graph  $\{y=\psi(x), x\in\Delta_{r}^{-}\}$ for some $r>0$ sufficiently small. Similarly, when $\lambda=e^{2\pi i p/q}$, there are $q$  fat repelling petals $\mathcal{P}_{{\rm rep}, j}$, $1\leq j\leq q$, which can be defined as the preimages under $T$ of a smooth graph $\{y=\psi(x), x\in \Delta_{r}^{-}\}$, for $r>0$ sufficiently small.

Let $\Patt$ and $\Prep$ denote the union of the $q$ attractive, and respectively of the $q$ repelling petals. 
We have $f(\overline{\Patt})\subset \Patt\cup W^{ss}(0)$ and all points in $\Patt$ are attracted to the semi-parabolic fixed point at the origin in forward time. Similarly $f^{-1}(\overline{\Prep})\subset \Prep\cup\{0\}$ and all points in $\Prep$ are attracted to the origin in backward time.

The repelling and attracting petals, sliced by a center manifold, alternate as we turn around the origin. By construction, the intersection of an attracting petal with the two adjacent repelling petals consists of two open disjoint disks. 
More precisely, let $D^{1}_{j}$ and $D^{2}_{j}$ be the two connected components of $\mathcal{P}_{{\rm att}, j}\cap \mathcal{P}_{\rm rep}$. Then 
\begin{equation}\label{eq:Pinv}
\mathcal{P}^{k}_{{\rm inv},j}=\bigcup_{n\in\Z} f^{nq}(D^{k}_{j}),\ \ \mbox{for}\ k=1, 2\ \ \mbox{and}\ \ 1\leq j\leq q
\end{equation}
are $2q$ completely invariant local petals. They are contained in the union of $\Prep$ and $\Patt$. We view these invariant petals as small, because they are {\it a~priori} defined only on the neighborhood $U$ where the map $f$ is conjugate to the normal form \eqref{eq:NF3}.

In Section \ref{sec:proof} we prove the existence of  hedgehogs for holomorphic germs of diffeomorphisms of $(\C^2,0)$ with a semi-neutral fixed point at the origin. Following the original strategy of P\'erez-Marco, we construct the hedgehog as a Hausdorff limit of invariant petals for approximating germs with a semi-parabolic fixed point at $0$. However, the normalizing domains on which the semi-parabolic germs can be conjugate to their corresponding normal forms given in Proposition \ref{thm:normalform3} necessarily shrink to $0$ as the sequence converges to a semi-Cremer germ, because the semi-Cremer germ is non-linearizable. Therefore we first need to construct big invariant petals before applying the construction in Section \ref{sec:proof}. 

Let $B\Subset B'$ be as in Theorem \ref{thm:Hedgehog} such that $f$ is partially hyperbolic on $B'$. There exist a horizontal cone field $\mathcal{C}^h$ which is forward invariant and a vertical cone field $\mathcal{C}^v$  which is backward invariant on $B'$, as in Equation \eqref{eq:inv}. We  say that an analytic curve $\gamma$ is {\it vertical-like}/{\it horizontal-like} if for any point $y$ on $\gamma$, the tangent space to $\gamma$ at $y$ is contained in the vertical/horizontal cone at $y$. Let $W^{ss}_{\rm loc}(0)$ be the local stable manifold of the semi-indifferent fixed point, i.e. the connected component of $W^{ss}(0)\cap B$ which contains $0$. We further assume that $W^{ss}_{\rm loc}(0)$ is vertical-like.

We show that the invariant petals are big with respect to the ball $B$ (a natural way to express this is to ask that they touch the boundary of $B$). 
In dimension one, P\'erez-Marco achieves this using the Uniformization Theorem and analytic circle diffeomorphisms, tools which are not readily available in the two-dimensional setting. 
The main obstruction is that, while each local invariant petal is contained in a complex analytic line (the asymptotic curve $\Sigma$ introduced below), the union of the $2q$ petals belongs to a (non-unique) center manifold of class $C^k$ for some $k\geq 1$, which is not complex or real analytic (see e.g. \cite{vS}). Instead of complex methods we will use some topological tools: Brouwer's Plane Translation Theorem \ref{thm:Brouwer} and covering space theory.

Using Theorem \ref{thm:InvPetals}, we define the asymptotic curve(s) $\Sigma$ to be the set of points in the domain of $f$, different from $0$, which are attracted to $0$ under backward iterations of $f$. The set $\Sigma$ has $q$ connected components which contain $0$ in the boundary, and each is a $f^{q}$-invariant Riemann surface immersed in $\C^2$. If $f$ is a global diffeomorphism of $\C^2$, Ueda \cite{U2} showed that these are biholomorphic to $\C$.

\begin{prop}\label{prop:horizontal} The set $\Sigma_{B}=\{x\in \Sigma : f^{-n}(x)\in B\ \mbox{for all}\ n\geq 0\}$ is horizontal-like.
\end{prop}
\proof
Let $x\in\Sigma_{B}$. There exists $m>0$ such that $y=f^{-m}(x)\in \Prep$. The repelling petal $\Prep$ is horizontal-like from \cite{U2} and the construction above. Therefore any tangent vector to $\Sigma$ at $y$ belongs to the horizontal cone at $y$, and all forward iterates $f^{i}(y), 0\leq i\leq m$, remain in $B$, so $T_{x}\Sigma$ is contained in the horizontal cone at $x$.
\qed

\begin{figure}[htb]
\begin{center}
\includegraphics[scale=0.277]{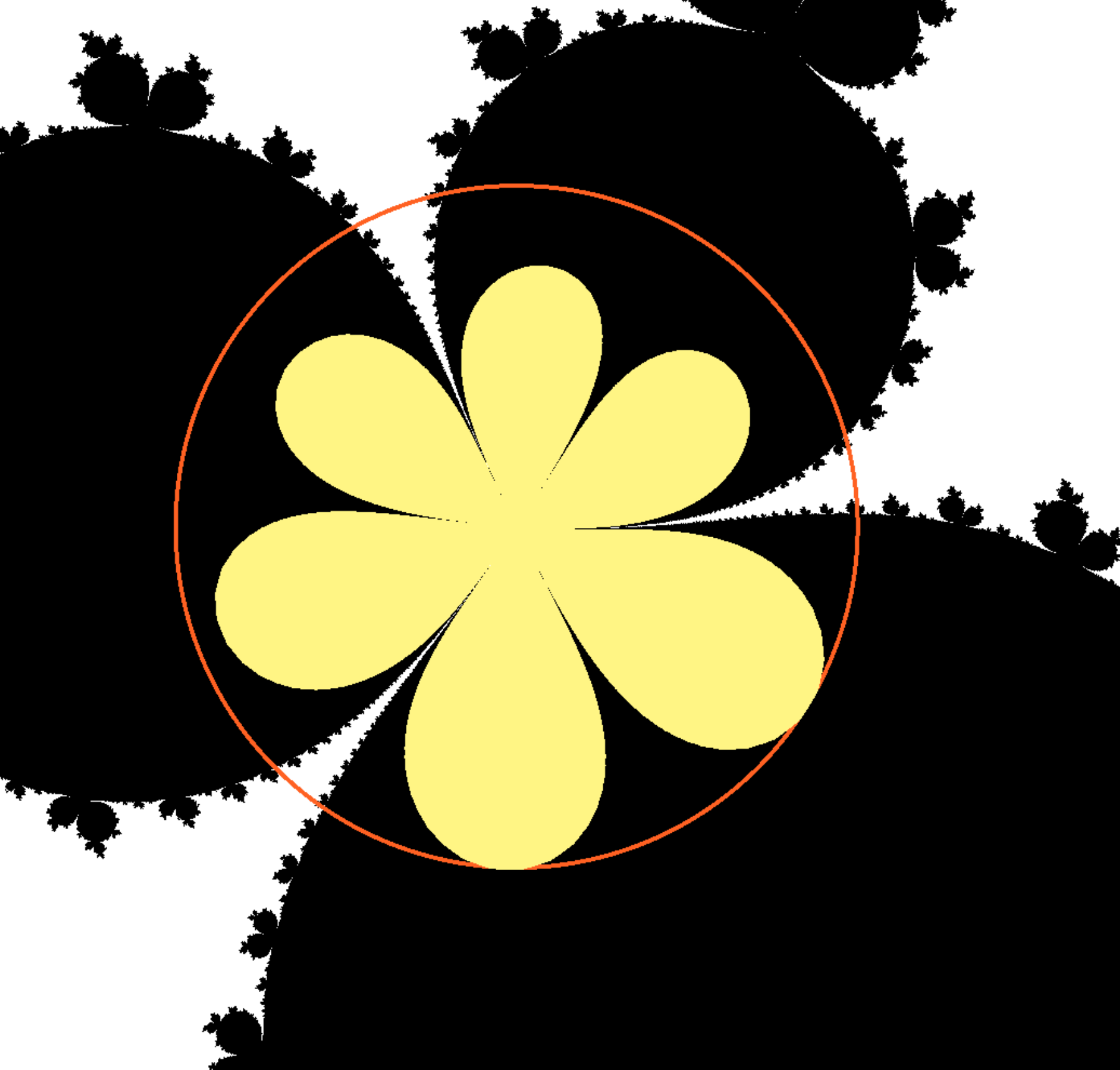}
\end{center}
\caption{Maximal invariant petals for $q=3$, relative to the ball $B$. Some petals touch the boundary of $B$.}
\label{fig:petals}
\end{figure}

Let $\mathcal{P}$ be the union of the $2q$ connected components of the set
\begin{equation}\label{eq:Pmax}
\{x\in B\setminus \{0\} : f^{ n}(x)\in B\ \forall n\in\Z\ \ \mbox{and}\ f^{n}(x)\rightarrow0\ \mbox{as}\ n\rightarrow \pm\infty \}
\end{equation}
which contain $0$ in their boundaries. We refer to $\mathcal{P}$ as the set of maximal invariant petals relative to the ball $B$. 
By definition $\mathcal{P}\subset \Sigma_B$, hence it is horizontal-like by Proposition \ref{prop:horizontal}. 
Each component of $\mathcal{P}$ contains a local invariant petal as defined by Equation \eqref{eq:Pinv} and Theorem \ref{thm:InvPetals} and is invariant by $f^{q}$ (see Figure \ref{fig:petals}).

Let $\overline{\mathcal{P}}$ and $\partial \mathcal{P}$ denote the closure, respectively the boundary of the set $\mathcal{P}$ in $\C^{2}$. 
In the following two propositions we collect a couple of elementary results about $\mathcal{P}$.

\begin{prop}\label{prop:Popen} 
The set $\mathcal{P}$ is open rel $\Sigma$ and its connected components are simply connected.
\end{prop}
\proof 
Let us first notice that the only point of intersection of the vertical-like local strong stable manifold $W^{ss}_{\rm loc}(0)$ and the horizontal-like set $\overline{\mathcal{P}}$ is $0$, by transversality. If $x\in W^{ss}(0)\cap \overline{\mathcal{P}}$ then there exists a positive integer $m$ such that $f^{m}(x)\in W^{ss}_{\rm loc}(0)\cap \overline{\mathcal{P}}$, which implies that $x$ is the fixed point $0$, which does not belong to $\mathcal{P}$. Therefore $W^{ss}(0)\cap \mathcal{P}=\emptyset$.

We now show that $\mathcal{P}$ is open relative to $\Sigma$. Let $\mathcal{B}_{\rm par}(0)= \bigcup_{n\geq 0}f^{-n}(\Patt)$ be the basin of the semi-parabolic fixed point $0$.  If $x\in \mathcal{P}$ then $x\in \mathcal{B}_{\rm par}(0)$, which is open in $\C^{2}$. Moreover, since $\Patt$ and  $\Prep$ are bases of convergence for $f$ on $\mathcal{B}_{\rm par}(0)$ and respectively for $f^{-1}$ on $\Sigma$, there exists a first iterate $n$ such that $f^{n}(x)\in \Patt$ and a first iterate $m$ such that $f^{-m}(x) \in \Prep$. There exists a neighborhood $U\subset\Sigma$ of $x$ such that $f^{n}(U)\subset \Patt$, $f^{-m}(U)\subset \Prep$ and $\bigcap_{-m<i<n}f^{-i}(U)\subset B$. Hence $U$ is an open set (rel $\Sigma$) contained in $\mathcal{P}$.

An immediate application of the Maximum Modulus Principle shows that each connected component of the set $\mathcal{P}$ is simply connected. Suppose that $\gamma$ is a non-trivial loop in $\mathcal{P}$. The set $\mathcal{P}$ is contained in $\Sigma$, so by eventually considering an iterate $f^{-nq}(\gamma)$ we may assume that $\gamma$ is contained in a local repelling petal $\mathcal{P}_{{\rm rep}, j}$ which is simply connected. Let $D$ denote the small disk bounded by $\gamma$ in $\mathcal{P}_{{\rm rep}, j}$. For each $n$, the map $\|f^{n}\|^{2}$ is subharmonic on $D$, so it attains its maximum on $\gamma$.
This implies that for every point in $D$, all forward and backward iterates belong to $B$, so $\{f^{n}\}_{n\in \Z}$ is a normal family on the open set $\mathcal{P}\cup D$. Let $f^*$ be the limit of any convergent subsequence. $f^*$ is a holomorphic function which is identically $0$ on $\mathcal{P}$, hence it must vanish identically on $D$ as well. In conclusion, for every point in $D$, all forward and backward iterates converge to $0$. It follows that $D$ belongs to $\mathcal{P}$, hence $\gamma$ is null homotopic and henceforth $\mathcal{P}$ is simply connected. 
\qed

\begin{prop}\label{prop:E}~\
\begin{itemize}
\item[a)] $\overline{\mathcal{P}}$ and $\partial \mathcal{P}$ are completely invariant.
\item[b)] Let $W^c(0)$ be any center manifold of the semi-parabolic fixed point, defined locally around $0$ as graph of a function $\varphi_f: E^c \cap B' \rightarrow E^s$. The set $\overline{\mathcal{P}}$ is connected, and compactly contained in $W^c(0)$.
\item[c)] The boundary $\partial \mathcal{P}$ cannot contain any attracting or hyperbolic fixed points of $f^q$.
\end{itemize}
\end{prop}
\proof
By construction, the set $\mathcal{P}$ is completely invariant, hence its closure and its boundary are also completely invariant.
By the weak uniqueness property of center manifolds, the set $\overline{\mathcal{P}}$ belongs to every center manifold $W^c(0)$ locally defined as a graph of a function on $E^c\cap B'$. 

Since we do not know the dynamics of $f$ on the boundary of $\mathcal{P}$, we cannot {\it a priori} assume that the boundary of $\mathcal{P}$ in contained in $\Sigma\cup\{0\}$. In any case, the sets $\partial{\mathcal{P}}$ and $\overline{\mathcal{P}}$ are contained in $\overline{B}$ and 
compactly contained in $W^c(0)$, since $W^{c}(0)$ is properly embedded into the bigger ball $B'$. The connectedness of $\overline{\mathcal{P}}$ follows from the fact that each connected component of $\mathcal{P}$ contains $0$ in its boundary.

For part c), assume that $z\in\partial\mathcal{P}$ is a fixed point of $f^{q}$ different from $0$. If $z$ is attracting then its basin of attraction $\mathcal{B}_{\rm att}(z)$ is open and points in $\mathcal{B}_{\rm att}(z) \cap \mathcal{P}$  converge both to $z$ and to $0$ under forward iterations of $f^q$, which is impossible. If $z$ is hyperbolic, then it has a stable and an unstable manifold. Let 
\begin{equation*}\label{eq:Wuloc}
W^{u}_{\rm loc}(z)=\{x\in B': f^{-nq}(x) \in B'\ \forall n\in\N\ \ \mbox{and}\ f^{-nq}(x)\rightarrow z \mbox{ as } n\rightarrow\infty\}
\end{equation*} 
denote the local unstable manifold relative to $B'$. Since $z$ belongs to the  center manifold(s) at $0$, by the weak uniqueness property, $W^{u}_{\rm loc}(z)$ belongs to all $W^c(0)$, so they must coincide in a neighborhood of $z$. This is again a contradiction, since any neighborhood of $z$ in $W^c(0)$ contain points from $\mathcal{P}$, which will therefore be forced to converge under backward iterations to $z$, as well as to $0$. This is impossible, since $z\neq 0$. 
\qed

We can use the local dynamics as in the proof of Proposition \ref{prop:E} to exclude the possibility of having a semi-Siegel or another semi-parabolic fixed point in the boundary of $\mathcal{P}$. The only case that one cannot elementary exclude, is the existence of a semi-Cremer fixed point in $\partial\mathcal{P}$.

We now show that the closure of $\mathcal{P}$ meets the boundary of the ball $B$. For this we need some topological tools about homeomorphisms of the plane. Consider a fixed point free orientation-preserving homeomorphism $h$ of $\R^2$.

\begin{defn}[\textbf{Domain of translation}]\label{def:DT} A domain of translation for $h$ is an open connected subset of $\R^{2}$ whose boundary is $\ell\cup h(\ell)$ where $\ell$ is the image of a proper embedding of $\R$ in $\R^{2}$, such that $\ell$ separates $h^{-1}(\ell)$ and $h(\ell)$. 
\end{defn}

Let $ U(D)=\bigcup_{n\in\Z}h^n(D)$, where $D$ be the closure of a domain of translation for $h$. The set $U(D)$ is open and connected and $h:U(D)\rightarrow U(D)$ is conjugate to the translation $T:\R^2\rightarrow \R^2$ given by $T(x,y)=(x+1,y)$.

\begin{thm}[\textbf{Brouwer's Plane Translation Theorem \cite{Fr}}]\label{thm:Brouwer} Suppose that $h:\R^{2}\rightarrow \R^{2}$ is an orientation-preserving homeomorphism of the plane without fixed points. Then every point is contained in some domain of translation.
\end{thm}

\begin{lemma}\label{lemma:boundary-pt}
Let $U\subset \R^{2}$ be a simply connected bounded domain. Let $h$ be an orientation-preserving homeomorphism of $\R^{2}$ without fixed points in $U$.
Suppose $C\subset U$ is a connected, invariant  set for $h$ such that $\overline{C} \cap\partial U =\{z_{0}\}$. Then $z_0$ is a fixed point of $h$ and $h^{n}(z)\rightarrow z_{0}$ as $n\rightarrow \pm \infty$, for all $z\in C$.
\end{lemma}
\proof
Let $O$ be the connected component of $\R^{2}\setminus {\rm Fix}(h)$ that contains $U$, where ${\rm Fix}(h)$ denotes the fixed points of the homeomorphism $h$. Since $C$ is invariant by $h$, we have $h(O)=O$. Let $\tilde{O}$ be the universal cover of $O$ and $p:\tilde{O}\rightarrow O$ be the covering map. Then $\tilde{O}$ is homeomorphic to $\R^{2}$. In fact, the only simply connected surfaces without boundary are $\s^{2}$ and $\R^{2}$, but $\s^{2}$ cannot be the universal cover of $O$ since $O$ is not compact.  

Let $x_{0}$ be a point in $C$ and let $\tilde{x}_{0}$ be a point in the fiber $p^{-1}(x_{0})$. Let $i:U\hookrightarrow O$ be the inclusion of $U$ into $O$. By hypothesis, the set $U$ is open and simply connected, hence there exists a lift $s:(U,x_{0})\rightarrow (\tilde{O},\tilde{x}_{0})$ such that $p\circ s = i$. Let $\tilde{U}$ be the connected component of $p^{-1}(U)$ which contains $\tilde{x}_{0}$. Then $\tilde{U}=s(U)$ and the restriction $p:\tilde{U}\rightarrow U$ is a homeomorphism whose inverse is the section $s$. 

We  now lift the homeomorphism $h:O\rightarrow O$ to the universal cover. The points $x_{0}$ and $h(x_{0})$ belong to $C$, so they also belong to $U$. 
Let $\tilde{x}_{1}=s(h(x_{0}))$. 
By general covering space theory, there exists a lift $\tilde{h}:\tilde{O}\rightarrow \tilde{O}$ such that $p \circ \tilde{h}=h\circ p$ and $\tilde{h}(\tilde{x}_{0})=\tilde{x}_{1}$. 

The map $\tilde{h}$ is a homeomorphism of the plane $\tilde{O}$ because it is a lift which induces a homeomorphism on fundamental groups. Moreover, $\tilde{h}$ is orientation-preserving and does not have any fixed points since $h$  is orientation-preserving and does not have any fixed points in $O$.  
By Brouwer's Plane Translation Theorem \ref{thm:Brouwer}, every point in $\tilde{O}$ is contained in some domain of translation. 

Suppose $\overline{C} \cap\partial U =\{z_{0}\}$. Let $\tilde{C} = s(C)$. Then $\tilde{h}(\tilde{C})$ and $\tilde{C}$ are two non-disjoint connected components of $p^{-1}(C)$, so they must coincide. Therefore $\tilde{C}$ is invariant by $\tilde{h}$. Let $z$ be any point in $C$. Then $\tilde{z}=s(z)$ is a point in $\tilde{C}$. Let $D$ be the closure of  a domain of translation containing $\tilde{z}$. 
The restriction of $\tilde{h}$ to $U(D)= \bigcup_{n\in \Z}\tilde{h}^{n}(D)$ is topologically conjugate to the translation  $T$ on $\R^{2}$. This shows that $\tilde{h}^{n}(\tilde{z})\rightarrow\infty$ as $n\rightarrow\pm\infty$. The homeomorphism $h$ on $C$ is conjugate to $\tilde{h}$ on $\tilde{C}$. The sequence of iterates $h^{n}(z)$ has to converge to a point on the boundary of the set $O$, but the only point on the boundary of $O$ and $C$ is $z_{0}$. So $h^{n}(z)\rightarrow z_{0}$ as $n\rightarrow \pm\infty$. 
\qed

\begin{remark}\label{rem:CL}
The idea of working in a covering space of a component of the complement of the set of fixed points,
as in the proof of Lemma \ref{lemma:boundary-pt}, 
was also used by Brown \cite{B} to reprove the Cartwright-Littlewood Theorem.
\end{remark}

\begin{thm}[\textbf{Maximal petals relative to $B$}]\label{thm:PBd} If the boundary $\partial \mathcal{P}$ does not contain a fixed point of $f^{q}$ other than the semi-parabolic fixed point $0$, then $\overline{\mathcal{P}}\cap \partial{B}\neq \emptyset$.
\end{thm}
\proof
We show that if $\partial \mathcal{P}$ does not contain any fixed point of $f^{q}$ different from $0$, then every point of the boundary converges to $0$ under forward and backward iterations. 

\begin{figure}[htb]
\begin{center}
\includegraphics[scale=0.9]{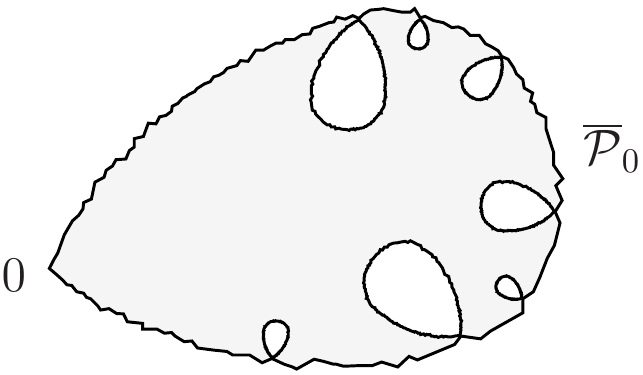}
\end{center}
\caption{An illustration of a set $\mathcal{P}_{0}$ with a complicated boundary.}
\label{fig:intricate}
\end{figure}

By modifying $f$ if needed with a cut-off function outside $B'$, as discussed in Section \ref{sec:centermfd}, we can consider a global center manifold $W^{c}(0)\simeq \R^{2}$. 
Let $\mathcal{P}_0$ be any connected component of $\mathcal{P}$. It is invariant by $f^q$ and by Proposition \ref{prop:E}  $\overline{P}_0$ is compactly contained in $W^c(0)$. By Proposition \ref{prop:Popen}, $\mathcal{P}_{0}$ is open rel $\Sigma$ (hence also open rel $W^c(0)$) and simply connected. The map $f^q$ is analytic on $\mathcal{P}_{0}$, hence orientation-preserving.
The boundary of $\mathcal{P}_{0}$ contains no fixed points of $f^{q}$ different from $0$ by assumption, but it could presumably be a complicated topological set (see Figure \ref{fig:intricate}). It could be non-locally connected or even the common boundary of three disjoint connected open sets (lakes of Wada). 

In what follows, set $h:=f^{q}$. Let $\Omega=\R^{2}\cup\{\infty\}\setminus \overline{\mathcal{P}}_{0}$ and denote by $\Omega_{0}$ the unbounded connected component of $\Omega$. The set $\Omega$ is an open subset of the sphere. By a standard result in topology, the set $\R^{2}\cup\{\infty\}\setminus\Omega$ is connected if and only if each connected component of $\Omega$ is simply connected. The closure of $\mathcal{P}_0$ is connected, hence $\Omega_{0}$ is simply connected, which implies that $\partial{\Omega_{0}}$ is a connected subset of $\partial{\mathcal{P}}_{0}$. Let $C = \partial{\Omega}_{0}$ denote the outer boundary of the set $\mathcal{P}_{0}$. It is compact, connected, and invariant under $h$. To simplify notation, let $C^{*}$ denote $C\setminus \{0\}$.

\begin{remark}\label{rem:analytic}
In principle, $C$ is the entire boundary of $\mathcal{P}_{0}$ (this would certainly be true is $h$ were analytic on $\overline{\mathcal{P}}_{0}$, by normality and the maximum modulus principle), but since the center manifold is not analytic, we do not know the topology of $\overline{\mathcal{P}}_{0}$, 
so we have to assume the most general situation. 
\end{remark}

\begin{figure}[htb]
\begin{center}
\includegraphics[scale=0.87]{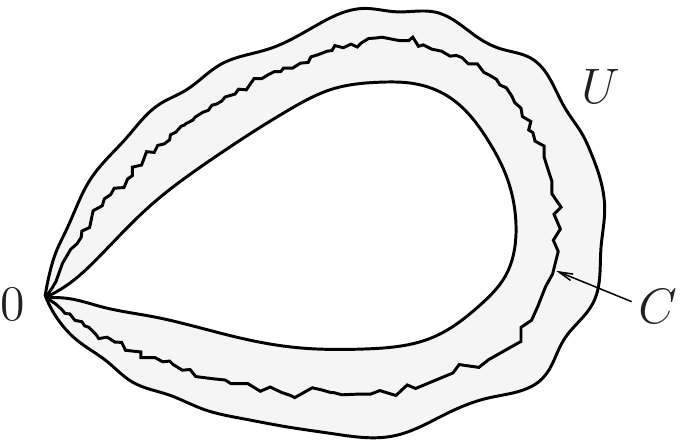}
\end{center}
\caption{A croissant $U$ containing $C^{*}$, where $C$ is the outer boundary of the invariant petal $\mathcal{P}_{0}$.}
\label{fig:nbdU}
\end{figure}

From the dynamics of the semi-parabolic germ, we know that there is a neighborhood around $0$ on which $h(x)\neq x$ whenever $x\neq 0$. By assumption $h$ does not have any fixed points on $C^{*}$. By continuity of the map $h$, there exists an open, connected and simply connected set $U\supset C^{*}$ in the shape of a croissant, as shown in Figure \ref{fig:nbdU}, such that $h$ does not have any fixed points on $U$. The set $C$ intersects the boundary of $U$ only at $0$.  The boundary of $U$ consists of two simple closed curves $C_{1}$ and $C_{2}$ touching at $0$.

Since $C_{1}$ and $C_{2}$ are simple closed curves, the Jordan-Schoenflies Theorem (see \cite{C}) allows us to extend the homeomorphism $h$ to a homeomorphism of the whole plane (which we also denote $h$). However, the homeomorphism $h$ constructed like this may have other fixed points outside of the domain $U$.

We now apply Lemma \ref{lemma:boundary-pt} for the homeomorphism $h$ and the invariant set $C^{*}\subset U$. Hence $h^{n}(z)\rightarrow 0$ as $n\rightarrow \pm\infty$, for any point $z\in C^{*}$. 

Suppose $\overline{\mathcal{P}}\cap \partial{B}= \emptyset$. Then $\overline{\mathcal{P}} \subset \{x\in B : f^{ n}(x)\in B\ \forall n\in\Z\}$.
However, we have shown that all points on  the outer boundary of $\mathcal{P}$, converge to $0$ under forward and backward iterations. So $\partial \mathcal{P}\subset \mathcal{P}$, which is open (rel $\Sigma$).
This is a contradiction, so $\overline{\mathcal{P}}$ must intersect the boundary of the ball $B$. 
\qed

From the proof of Theorem \ref{thm:PBd} and the local dynamics around the semi-parabolic fixed point we also have the following immediate consequence:
\begin{prop}\label{prop:closure} 
The set $\overline{\mathcal{P}}$ is simply connected.
\end{prop}

\section{Proof of the main theorem}\label{sec:proof}

We now have all the ingredients to complete the proof of Theorem \ref{thm:Hedgehog}.  

\begin{lemma}\label{lemma:approx} Let $f$ be a  germ of holomorphic diffeomorphisms of $(\C^{2},0)$ with a semi-indifferent fixed point at $0$ with eigenvalues $\lambda=e^{2\pi i \alpha}$, $\alpha\notin \Q$, and $\mu$, with $|\mu|<1$. Suppose $p_n/q_n$ are the convergents of $\alpha$ given by the continued fraction algorithm. Consider a neighborhood $B'\subset \C^{2}$ of $0$, on which $f$ is partially hyperbolic. There exists a sequence $f_{n}\rightarrow f$ of germs of holomorphic diffeomorphisms of $(\C^{2},0)$ such that
\begin{itemize}
\item[a)] $f_{n}$ has a semi-parabolic fixed point at $0$, with eigenvalues $\lambda_n$ and $\mu_n$,  $\lambda_{n}=e^{2\pi i p_{n}/q_{n}}$ and $|\mu_n|<1$, of semi-parabolic multiplicity $1$.
\item[b)] $f_{n}$ does not have other semi-neutral periodic points of period $\leq q_{n}$ inside $B'$.  
\item[c)] $f_{n}$ is partially hyperbolic on $B'$ for $n$ large enough.
\end{itemize} 
\end{lemma}
\proof 
Consider a sequence of semi-parabolic germs 
\[
f_{n}(x,y)=\left(\lambda_n x+ h.o.t., \mu_{n} y + h.o.t.\right)
\] 
converging uniformly to $f$ on the neighborhood $B'$, with Jacobian matrix at 0 equal to $\mbox{Diag}(\lambda_n, \mu_{n})$, where $\lambda_{n}=e^{2\pi i p_{n}/q_{n}}$ and $0<|\mu_n|<1$.
If $f_{n}$ does not satisfy property b) for some $n$, then we can make an arbitrarily small change of the germ $f_{n}$ while keeping the eigenvalue $\lambda_n$ fixed, so that any semi-neutral periodic point (different from the origin) of period $\leq q_{n}$ inside $B'$ becomes attracting or hyperbolic. 
Denote for simplicity the new germ also by $f_{n}$. 

By perturbing the germ $f_{n}$ slightly if necessary while keeping the eigenvalue $\lambda_n$ fixed, we can ensure that the semi-parabolic multiplicity of the origin is one. In conclusion, there exists a sequence $f_{n}\rightarrow f$ satisfying the first two claims of the lemma. The convergence is uniform, so we can assume that  $f_{n}$ is partially hyperbolic on $B'$, for $n$ sufficiently large.
\qed

We are now able to establish the existence of the set $\mathcal{H}$ from the main theorem and obtain some immediate properties of this set.

\medskip
\noindent\textbf{Proof of Theorem \ref{thm:Hedgehog}.}  Let $f$ be a germ of holomorphic diffeomorphisms of $(\C^{2},0)$
with a semi-indifferent fixed point at $0$ with eigenvalues $\lambda$ and $\mu$, where $|\lambda|=1$ and $|\mu|<1$.  Let $f_{n}$ be a sequence of semi-parabolic germs converging uniformly to  $f$ as in Lemma \ref{lemma:approx}. 

Consider a ball $B\subset \C^2$ centered at $0$ with $\overline{B}\subset B'$ and denote by $\mathcal{P}_{n}$ the maximal invariant petals of the semi-parabolic germ $f_{n}$ relative to the ball $B$ (see Equation \eqref{eq:Pmax}). Let $ \mathcal{H}$ be the limit of a convergent subsequence $(\overline{\mathcal{P}}_{n_k})_k$ in the Hausdorff topology of compact subsets of $\C^{2}$. Since each set $\overline{\mathcal{P}}_{n}$ is compact, connected and completely invariant under $f_n$, the Hausdorff limit of any convergent subsequence of $(\overline{\mathcal{P}}_{n})_n$ will also be compact, connected, and completely invariant under $f$. 

By part b) of Lemma \ref{lemma:approx} and part c) of Lemma \ref{prop:E} it follows that $\partial \mathcal{P}_n$ cannot contain any fixed points of $f^{q_n}$ other than $0$. Therefore, we can apply Theorem \ref{thm:PBd} to conclude that $\overline{\mathcal{P}}_{n}\cap \partial B\neq \emptyset$ for every $n$, hence $\mathcal{H}\cap\partial B \neq \emptyset$. 

Part a) of the theorem is an easy consequence of the weak uniqueness property of center manifolds, since all center manifolds defined relative to the neighborhood $B' \supset \overline{B}$ contain the maximal invariant set of $f$ in $B$. The existence of the strong stable foliation in part d) follows by general theory of partially hyperbolic systems.

Let $\tilde{f}$ be a smooth extension of $f$ to $\C^{2}$ such that $\tilde{f}=f$ on $B'$, as in Section \ref{sec:centermfd}. Let $W^c(0)$ be the global center manifold of $\tilde{f}$. 
 If $\mathcal{H}$ is not full relative to $W^c(0)$, then we will replace it by $\hat{\mathcal{H}}=\mathcal{H}\cup \Omega$, where $\Omega$ is the union of the bounded connected components of the complement of $\mathcal{H}$ in $W^c(0)$. Clearly $\tilde{f}(\Omega)=\Omega\subset B'$, which implies that $f(\Omega)=\Omega$. Since $f$ is a diffeomorphism on $B'$, it follows that $f^{-n}(\Omega)=\Omega$ for all $n\geq 1$. In particular all backward iterates of $\Omega$ never escape the set $B'$, hence using the weak uniqueness property once more, we conclude that the set $\Omega$ belongs to all center manifolds defined relative to $B'$. Therefore, the set $\hat{\mathcal{H}}$ is full and satisfies all the other properties of the theorem. 
\qed

Note that the tools outlined in Section \ref{sec:semiparabolic} and in the current section also apply to one-dimensional germs of holomorphic diffeomorphisms of $(\C,0)$ with a neutral fixed point at $0$. In this context, the role of the center manifold is taken by the ambient complex plane $\C$. Using Theorem \ref{thm:PBd} and Lemma \ref{lemma:boundary-pt}, we can construct the hedgehog $\mathcal{H}$ relative to a domain $B$ as a Hausdorff limit of a sequence of maximal invariant petals defined relative to $B$, corresponding to a sequence of parabolic germs $f_n$ converging to $f$. This gives an alternative proof to the theorem of P\'erez-Marco, stated below.

\begin{thm}[P\'erez-Marco \cite{PM1}]\label{thm:PerezMarco} Let $f(z)=\lambda z + \bigO(z^{2})$, with $|\lambda|=1$, be a local holomorphic diffeomorphism, and $B$ a Jordan domain around the neutral fixed point $0$. Assume that $f$ and $f^{-1}$ are defined and univalent in a neighborhood of $\overline{B}$. There exists a compact, connected, full set $\mathcal{H}$ containing $0$, such that $\mathcal{H}\cap \partial B \neq \emptyset$, and which is completely invariant under $f$.
\end{thm} 

\noindent{\bf Associated circle homeomorphism.} We can associate to $(f,\mathcal{H})$ a homeomorphism of the unit circle as in \cite{PM1}. 
Consider a global center manifold $W^c(0)$ which contains the hedgehog $\mathcal{H}$. 
$W^c(0)$ is homeomorphic to 
$\R^2$, and we can put a reference complex structure on the global center manifold and identify it with $\C$. The map $f$ induced on the copy of $\C$ will only be $C^1$-smooth. By the Uniformization Theorem, there exists a Riemann map $\psi:\hat{\C}-\overline{\D}\rightarrow \hat{\C}-\mathcal{H}$ with $\psi(\infty)=\infty$. The map $g = \psi^{-1}\circ f \circ \psi$ is a $C^1$ diffeomorphism in a neighborhood of $\mathbb{S}^1$ in $\hat{\C}-\overline{\D}$. Since $g$ is uniformly continuous in this neighborhood, it extends to a homeomorphism on $\mathbb{S}^1$ with rotation number $\rho_g$. To show that $\rho_g$ is equal to $\alpha$, let $p_n/q_n$ be the convergents of $\alpha$ given by the continued fraction algorithm. Consider a sequence of semi-parabolic germs $f_n$ converging to $f$ as in Lemma \ref{lemma:approx}. We can choose a family of $C^1$-smooth center manifolds $W^c(0)$ which depend $C^1$ on $f_n$. We build maximal invariants petals $\mathcal{P}_n$ relative to a ball $B$ of fixed size.  After eventually passing to a convergent subsequence, the sets $\mathcal{P}_n$ converge in the Hausdorff topology to $\mathcal{H}$, so the uniformizing maps $\psi_n:\hat{\C}-\overline{\D}\rightarrow \hat{\C}-\mathcal{P}_n$ converge to $\psi$ in the Carath\'eodory kernel topology. The corresponding circle homeomorphisms $g_n$ converge to $g$. The rotation number depends continuously on the function. It is easy to see that the rotation number of the circle homeomorphism $g_n$ associated to the semi-parabolic germ $f_n$ is $p_n/q_n$. Therefore, the rotation number of $g$ is the limit of $p_n/q_n$, so it is equal to $\alpha$. 

\medskip
The sets $\mathcal{P}_{n}$ used in the proof of Theorem \ref{thm:Hedgehog} are locally connected. However one cannot expect $\mathcal{H}$ to be locally connected, since this property is not preserved under taking Hausdorff limits. The following proposition supports this claim. 

\begin{prop}\label{prop:nlc} 
Let $f$ be a holomorphic germ of diffeomorphisms of $(\C^2,0)$ with an irrational semi-indifferent fixed point at the origin. Let $\mathcal{H}$ be the hedgehog constructed in Theorem \ref{thm:Hedgehog} and denote by ${\rm int^c}(\mathcal{H})$ the interior of $\mathcal{H}$ relative to a center manifold. 
If $0\notin{\rm int^c}(\mathcal{H})$, then $\mathcal{H}$ is not locally connected.
\end{prop}
\proof
Suppose that $\mathcal{H}$ is locally connected and consider a center manifold $W^{c}(0)$ of the fixed point $0$ which contains $\mathcal{H}$. We can identify $W^c(0)$ with $\C$ by a homeomorphism $\phi:W^c(0)\rightarrow \C$,  and denote by $K=\phi(\mathcal{H})$ the hedgehog in the new coordinates. By the Carath\'eodory Theorem, since $K$ is locally connected, the Riemann map $\psi :\C-\D\rightarrow \C-K$ has a continuous extension $\psi:\s^1\rightarrow K$. The function $f$ on the hedgehog $\mathcal{H}$ is therefore conjugate to an orientation-preserving homeomorphism $h:\s^1\rightarrow \s^1$ of the unit circle.  The map $h$ has a well defined rotation number $\alpha\notin \Q$, where $\alpha$ is the argument of the neutral eigenvalue of $Df_0$.  However, this is not possible, since by construction the map $h$ has a proper completely invariant closed set on $\s^1$, given by $\psi^{-1}\circ\phi(0)$.
\qed

\section{Appendix: Alternative approach in dimension one}\label{sec:App}

This section is of independent interest. 
We keep the same notations as in Section \ref{sec:semiparabolic}.  We present a more direct proof that the closure of the maximal invariant set $\mathcal{P}$ defined in \eqref{eq:Pmax} meets the boundary of the ball $B$, in the case when we know the topology of $\overline{\mathcal{P}}$ (that it is simply connected). This follows from the topological theory of parabolic germs in the plane of Le Roux \cite{LR}.

As observed in Remark \ref{rem:analytic}, a specific case when we already know that $\overline{\mathcal{P}}$ is simply connected is when $f$ is analytic on a neighborhood of $\overline{\mathcal{P}}$. This is true in dimension one. Let $f(z) = \lambda z +\bigO(z^{2})$ be a germ of holomorphic diffeomorphisms of $(\C,0)$ with a parabolic fixed point at $0$ of multiplier $\lambda=e^{2\pi i p/q}$ and consider $B$ a domain containing the origin such that $f$ and $f^{-1}$ are defined and univalent in a neighborhood of $\overline{B}$.
Thus in this section we give yet another way of showing that the local invariant petals of the one-dimensional parabolic germ $f$ extend to the boundary of $B$, which provides an alternative proof of Theorem \ref{thm:PerezMarco}.

Let $\mathcal{P}_{0}$ be a connected component of $\mathcal{P}$; it contains $0$ in its closure. As usual, we say that an open set $U\subset \R^{2}$ is a {\it Jordan domain} if its closure is homeomorphic to the closed unit disk in the plane.

\begin{prop}\label{prop:U-nbd} If the boundary $\partial \mathcal{P}_{0}$ does not contain any fixed points of $f^{q}$ other than $0$, then there exists a Jordan domain $U\subset \R^{2}$ containing $\overline{\mathcal{P}}_{0}$ such that $f^{q}(x)\neq x$ for all $x\in U\setminus\{0\}$. 
\end{prop}

Proposition 2.2 from \cite{LR}, stated below as Lemma \ref{lemma:extension}, is a useful topological result that allows us to extend the germ $f^{q}$ of $(U, 0)$ to a homeomorphism $h$ of $\R^{2}$ with a unique fixed point at $0$. The homeomorphism $h$ extends to a homeomorphism $h:\s^2\rightarrow\s^2$ with only two fixed points, at $0$ and at $+\infty$. We then use Theorem \ref{thm:LR} below to show that $\overline{\mathcal{P}}_{0}$ is contained in a translation domain for the homeomorphism $h$.

\begin{lemma}[\textbf{Extension lemma} \cite{LR}]\label{lemma:extension} Let $U$ and $V$ be two closed Jordan domains in $\R^{2}$ containing $0$ and $h:U\rightarrow V$ a homeomorphism with a unique fixed point at $0$. Suppose $W\subset U$ is an open, connected set containing $0$ such that  $h(W)\subset U$. Then there is a homeomorphism $\widehat{h}:\R^{2}\rightarrow \R^{2}$ which coincides with $h$ on $W$ and has only one fixed point at $0$. 
\end{lemma}

Consider $U$ as in Proposition \ref{prop:U-nbd} and denote $f^{q}$ by $h$. Since $\mathcal{P}_{0}$ is invariant by $f^{q}$, there exists an open, connected set $W\supset \mathcal{P}_{0}$ such that $W$ and $f^{q}(W)$ are both subsets of $U$.

Let $S$ denote the point $0$ and $N$ denote $+\infty$ on the sphere $\s^{2}=\R^{2}\cup\{+\infty\}$.  With a small abuse of notation, let $h:\s^{2}\rightarrow\s^{2}$ be an orientation-preserving homeomorphism with only two fixed points $\{S,N\}$ given by Lemma \ref{lemma:extension}. We now introduce some terminology for the dynamics of the homeomorphism $h$ on the sphere and refer to \cite{LR} for more details. An {\it attractive petal based at $S$} is a closed topological disk $P_{S}\subset \s^{2}$ such that $S\in \partial P_{S}$, $N\notin P_{S}$ and $h(P_{S})\subset {\rm int}(P_{S})\cup\{S\}$. A {\it repelling petal based at $S$} is an attractive petal based at $S$ for $h^{-1}$.  Similarly, one can define attractive and repelling petals based at $N$. 

\begin{defn}\label{def:croissants}
An {\it attractive croissant for the dynamics $N$-$S$} for the homeomorphism $h$ is a closed topological disk $D\subset \s^{2}$ such that $N, S\in \partial D$, $h(D)\subset {\rm int}(D)\cup\{S,N\}$ and for any neighborhood $W_{N}$ of $N$ there exists an attractive petal $P_{S}$ based at $S$ such that $P_{S}\subset D$ and $D\setminus P_{S}\subset W_{N}$. A {\it repelling croissant for the dynamics $S$-$N$} for $h$ is an attractive croissant for the dynamics $N$-$S$ for $h^{-1}$. 
\end{defn}

\begin{thm}[Le Roux \cite{LR}]\label{thm:LR} Let $h:\s^{2}\rightarrow\s^{2}$ be an orientation-preserving homeomorphism of the sphere with only two fixed points $\{S,N\}$ such that  $Index(N)=1-q<1$.  

There exists $q$ {\it attracting croissants} for the dynamics $N$-$S$ and $q$ {\it repelling croissants} for the dynamics $S$-$N$. The attracting and repelling croissants are cyclically alternating on the sphere and they intersect only at $S$ and $N$.  
\end{thm}

\begin{figure}[htb]
\begin{center}
\includegraphics[scale=0.71]{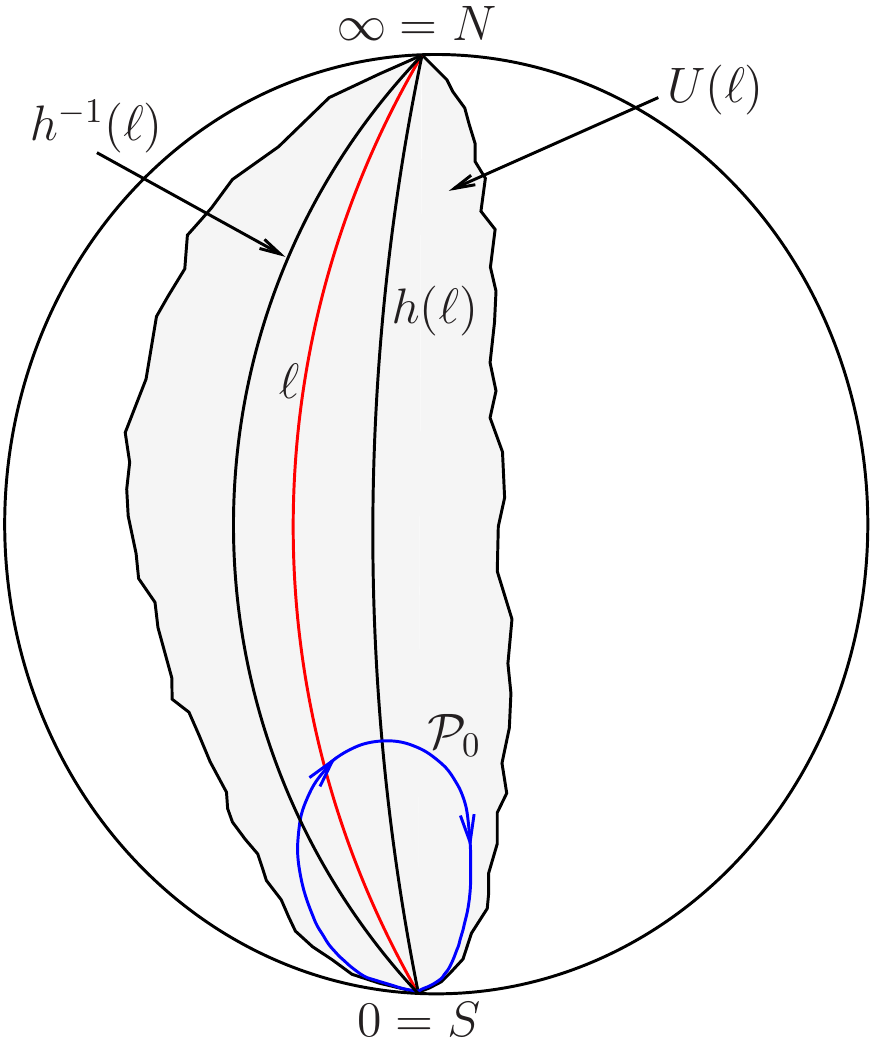}
\end{center}
\caption{The Brouwer domain $U(\ell)$ containing the invariant petal $\mathcal{P}_{0}$ on the sphere.}
\label{fig:sphere}
\end{figure}

The two lines that form the boundary of an attractive/repelling croissant and connect $S$ to $N$ are called {\it Brouwer lines} and they can be regarded as geodesics on the Euclidean sphere. Let $\ell$ be a Brouwer line bounding an attractive croissant for the dynamics $N$-$S$ given by Theorem \ref{thm:LR}.
So $\ell$ connects $S$ to $N$, but does not contain $S$ or $N$. Suppose that $h(\ell)$ is to the right of $\ell$ and $h^{-1}(\ell)$ is to the left of $\ell$. Let $D$ be the open topological disk on the sphere bounded by $\ell$ and $h(\ell)$. Define the domain $U(\ell)=\bigcup_{n\in\Z}h^{n}(D\cup\ell)$; it is called a Brouwer domain generated by $\ell$.  By  \cite[Section~3.2]{LR}, we know that there exists a homeomorphism 
$\phi:\R^{2}\rightarrow U(\ell)$ which conjugates the restriction $h:U(\ell)\rightarrow U(\ell)$ to the translation $T$ on $\R^{2}$, $T(x,y)=(x+1,y)$. So $U(\ell)$ is a domain of translation (see Figure \ref{fig:sphere}).

The index of the semi-parabolic fixed point is $1+q$ because we have not modified the homeomorphism in a small neighborhood around 0. By Lefschetz's formula, the sum of the indexes of the fixed points of $h$ is the Euler characteristic of the sphere, so the index at $S$ is $1+q$ if and only if  the index at $N$ is $1-q$.  So Theorem \ref{thm:LR} is applicable and in view of the discussion above, there exists some Brouwer line $\ell$ such that $\mathcal{P}_{0}\subset U(\ell)$. It follows that the semi-parabolic map $f^{q}$ is conjugated to a translation in a neighborhood of the invariant petal $\mathcal{P}_{0}$ inside the domain $U(\ell)$ on the sphere. Let $x\in \partial \mathcal{P}_{0}\setminus\{0\}$. This shows that there is a neighborhood $N_{\epsilon}(x)$ around $x$ such that if $y\in N_{\epsilon}(x)$ then $y$ converges to $0$ under forward or backwards iterations of $f^{q}$. In particular all points on the boundary of $\mathcal{P}_{0}$ converge to $0$ under forward or backwards iterations. 

This gives us a  mechanism of extending the invariant set $\mathcal{P}$ until it touches the boundary of $B$.



\begin{thebibliography}{99999}

\bibitem[B] {B} M. Brown, {\it A short short proof of the Cartwright-Littlewood theorem}, Proc. Amer. Math. Soc., 65 (1977), p. 372.

\bibitem[Brj] {Brj} A.D. Brjuno, Analytical form of differential equations. Transactions of the Moscow Mathematical Society {\bf 25}, 131-288 (1971); {\bf 26}, 199-239 (1972).

\bibitem[C] {C} S. Cairns, {\it An elementary proof of the Jordan-Schoenflies theorem}, Proc. Amer. Math. Soc., vol. 2 (1951), pp. 860-867.

\bibitem[Fr] {Fr} J. Franks, {\it A new proof of the Brouwer plane translation theorem}, Ergod. Th. and Dynam. Sys., 12 (1992), pp. 217-226.

\bibitem[Ha] {Ha}  M. Hakim, {\it Attracting domains for semi-attractive transformations of $\C^{p}$}, Publ. Mat. 38 (1994), no. 2, 479-499.

\bibitem[HPS] {HPS} M. Hirsch, C. Pugh, M. Shub, {\it Invariant manifolds}, Lecture Notes in Mathematics, vol. 583, Springer-Verlag, New York, 1977.

\bibitem[IS] {IS} H. Inou, M. Shishikura,  {\it The renormalization for parabolic fixed points and their perturbation}, Manuscript 2008.

\bibitem[LR] {LR} F. Le Roux, {\it Hom\'eomorphismes de surfaces - Th\'eor\`emes de la fleur de Leau-Fatou et de la vari\'et\'e stable}, Ast\'erisque 292 (2004). 

\bibitem[P] {P} Y. Pesin, {\it Lectures on Partial Hyperbolicity and Stable Ergodicity}. Z\"urich Lectures in Advanced Mathematics, EMS, 2004.

\bibitem[PM1] {PM1} R. P\'erez-Marco, {\it Fixed points and circle maps}, Acta Math., 179 (1997), 243-294.

\bibitem[PM2] {PM2} R. P\'erez-Marco, {\it On a question of Dulac and Fatou}, C. R. Acad. Sci. Paris, S\'erie I, 1995, 1045-1048.

\bibitem[PM3] {PM3} R. P\'erez-Marco, {\it Hedgehog's dynamics I}, Manuscript 2004.

\bibitem[PM4] {PM4} R. P\'erez-Marco, {\it Sur les dynamiques holomorphes non lin\'earisables et une conjecture de V. I. Arnol'd}, Ann. Sci. \'Ecole Norm. Sup. (4), 26(5):565-644, 1993.

\bibitem[R] {R} C. Robinson, {\it Dynamical Systems: Stability, Symbolic Dynamics, and Chaos}, 2nd ed., Studies in Advanced Mathematics, CRC Press, 1999.

\bibitem[RT]  {RT}  R. Radu, R. Tanase, {\it A structure theorem for semi-parabolic H\'enon maps}, 
{\tt arXiv:1411.3824 } 

\bibitem[R\"{u}s] {Rus} H. R\"{u}ssmann {\it \"{U}ber die Iteration analytischer Funktionen}, J Math. Mech. 17 (1967), pp. 523-532.

\bibitem[S] {S} M. Shub, {\it Global Stability of Dynamical Systems}, Springer-Verlag, 1987.

\bibitem[U1] {U1} T. Ueda, {\it Local structure of analytic transformations of two complex variables I}, J. Math. Kyoto Univ., 26(2) (1986), 233-261.

\bibitem[U2] {U2} T. Ueda, {\it Local structure of analytic transformations of two complex variables II}, J. Math. Kyoto Univ. 31 (1991), no. 3, 695-711.

\bibitem[V] {V}  A. Vanderbauwhede, {\it Centre manifolds, normal forms and elementary bifurcations}, in Dynamics Reported, A series in dynamical systems and their applications, Vol. 2, Wiley, Chichester, 1989, 89-169.

\bibitem[vS] {vS} S. van Strien, {\it Center manifolds are not $C^{\infty}$}, Math. Z. 166(2), 143-145 (1979).

\bibitem[Y1] {Y1} J. C. Yoccoz, {\it Analytic linearization of circle diffeomorphisms}, Dynamical systems and small divisors (Cetraro, 1998) Lecture Notes in Math., vol. 1784, Springer, Berlin, 2002, pp. 125-173.

\bibitem[Y2] {Y2} J. C. Yoccoz, {\it Petits diviseurs en dimension 1}, Ast\'erisque No. 231 (1995).
\end{thebibliography}
\end{document}